\overfullrule=0pt
\centerline {\bf Applying twice a minimax theorem}\par
\bigskip
\bigskip
\centerline {Biagio Ricceri}\par
\bigskip
\bigskip
\centerline {\it To Professor Wataru Takahashi, with esteem and friendship, on his 75th birthday}
\bigskip
\bigskip
{\bf Abstract:} Here is one of the results obtained in this paper:   Let $X, Y$ be two convex sets each in a real vector space, let $J:X\times Y\to {\bf R}$
be convex and without global minima in $X$ and concave in $Y$, and let $\Phi:X\to {\bf R}$ be strictly convex. Also, assume that, for some topology on $X$, $\Phi$ is lower semicontinuous
and, for each $y\in Y$ and $\lambda>0$, $J(\cdot,y)$ is lower semicontinuous and $J(\cdot,y)+\lambda\Phi(\cdot)$ is inf-compact.\par
Then,  for each $r\in ]\inf_X\Phi,\sup_X\Phi[$ and for each closed set $S\subseteq X$  satisfying
$$\Phi^{-1}(r)\subseteq S\subseteq \Phi^{-1}(]-\infty,r])\ ,$$
one has
$$\sup_Y\inf_SJ=\inf_S\sup_YJ\ .$$\par
\par
\bigskip
\bigskip
\bigskip
\bigskip
{\bf 1. Introduction}\par
\bigskip
A real-valued function $f$ on a topological space is said to be inf-compact (resp. sup-compact) if $f^{-1}(]-\infty,r])$ (resp. $f^{-1}([r,+\infty[$) is compact for all $r\in {\bf R}$.\par
\smallskip
A real-valued function $f$ on a convex set is said to be quasi-concave if 
$f^{-1}([r,+\infty[)$ is convex for all $r\in {\bf R}$.\par
\smallskip
In [3], we proved two general minimax theorems which, grouped together, can be stated as follows:\par
\medskip
THEOREM 1.A ([3], Theorems 1.1 and 1.2). - {\it Let $X$ be a topological space, $Y$ a convex set in a Hausdorff real topological vector space and
$f:X\times Y\to {\bf R}$ a function such that $f(\cdot,y)$ is lower semicontinuous, inf-compact and has a unique global minimum for all $y\in Y$.
Moreover, assume that either, for each $x\in X$, $f(x,\cdot)$ is continuous and quasi-concave or, for each $x\in X$, $f(x,\cdot)$
is concave.\par
Then, one has
$$\sup_Y\inf_Xf=\inf_X\sup_Yf\ .$$}\par
\medskip 
Theorem 1.A was first proved in the case where $Y$ is a real interval ([1], [2]) and successively extended to the general case
by means of a suitable inductive argument.\par
\smallskip
In [1], we applied Theorem 1.A (with $Y$ a real interval) to obtain a result ([1], Theorem 1) about the following problem: given two functions
$f, g:X\to {\bf R}$, find a interval $I\subseteq g(X)$ such that, for each $r\in I$, the restriction of $f$ to $g^{-1}(r)$ has a unique
global minimum.\par
\smallskip
The aim of the present paper is to establish a new minimax theorem (Theorem 2.1) which is the fruit of a joint application of Theorem 1.A and Theorem 1 of [1]. So, it follows, essentially, from a double application of Theorem 1.A, as the title stresses.\par
\smallskip
We then show some consequences of Theorem 2.1.\par
\bigskip
{\bf 2. Results}\par
\bigskip
In the sequel, $X$ is a  topological space, $Y$ is a non-empty set, $J:X\times Y\to {\bf R}$,
$\Phi:X\to {\bf R}$, $a, b$ are two numbers in $[0,+\infty]$,
with $a<b$. \par
\smallskip
For $y\in Y$ and $\lambda\in [0,+\infty]$,  we denote by $M_{\lambda,y}$ 
the set of all global minima of the function $J(\cdot,y)+\lambda\Phi(\cdot)$ if $\lambda<+\infty$, while
if $\lambda=+\infty$, $M_{\lambda,y}$  stands for the empty set.
We adopt the conventions $\inf\emptyset=+\infty$, $\sup\emptyset=-\infty$.
\smallskip
We also set
$$\alpha:=\sup_{y\in Y}\max\left \{ \inf_X \Phi,\sup_{M_{b,y}}\Phi\right \}\ ,$$
$$\beta:=\inf_{y\in Y}\min\left \{ \sup_X \Phi,\inf_{M_{a,y}}\Phi\right \}\ .$$
\smallskip
The following assumption will be adopted:\par
\smallskip
\noindent
$(a)$\hskip 5pt $Y$ is a convex set in a Hausdorff real topological vector space and either,
for each $x\in X$, the function $J(x,\cdot)$ is continuous and quasi-concave, or, for each $x\in X$, the function $J(x,\cdot)$ is concave.\par
\smallskip
Our main result is as follows:\par
\medskip
THEOREM 2.1. - {\it Besides $(a)$, assume that:\par
\noindent
$(a_1)$\hskip 5pt  $\alpha<\beta$\ ;\par
\noindent
$(a_2)$\hskip 5pt  $\Phi$ is lower semicontinuous\ ;\par
\noindent
$(a_3)$\hskip 5pt for each $\lambda\in ]a,b[$ and each $y\in Y$,  the function $J(\cdot,y)$ is lower semicontinuous
and the function $J(\cdot,y)+\lambda\Phi(\cdot)$ is inf-compact and admits a unique global minimum in $X$.\par
Then,  for each $r\in ]\alpha,\beta[$ and for each closed set $S\subseteq X$  satisfying
$$\Phi^{-1}(r)\subseteq S\subseteq \Phi^{-1}(]-\infty,r])\ ,\eqno{(2.1)}$$
one has
$$\sup_Y\inf_SJ=\inf_S\sup_YJ\ .\eqno{(2.2)}$$}\par
\smallskip
PROOF. Since $r\in ]\alpha,\beta[$, for each $y\in Y$, Theorem 1 of [1] (see Remark 1 of [1]) ensures the existence of $\lambda_{r,y}\in ]a,b[$ such that the unique
global minimum of $J(\cdot,y)+\lambda_{r,y}\Phi(\cdot)$, say $x_{r,y}$, lies in $\Phi^{-1}(r)$. Notice that $x_{r,y}$ is the only global minimum of the restriction of
the function $J(\cdot,y)$ to $\Phi^{-1}(]-\infty,r])$. Indeed, if not, there would exist $u\in \Phi^{-1}(]-\infty,r])$, with $u\neq x_{r,y}$, such that $J(u,y)\leq J(x_{r,y},y)$.
Then, (since $\lambda_{r,y}>0$) we would have
$$J(u,y)+\lambda_{r,y}\Phi(u)\leq J(x_{r,y},y)+\lambda_{r,y}\Phi(u)\leq J(x_{r,y},y)+\lambda_{r,y}r=J(x_{r,y},y)+\lambda_{r,y}\Phi(x_{r,y})$$
which is absurd. Therefore, since $S$ satisfies $(2.1)$, the restriction of $J(\cdot,y)$ to $S$ has a unique global minimum.
Now, observe that, for each $y\in Y$, $\rho\in {\bf R}$, $\lambda\in ]a,b[$, one has
$$\{x\in S : J(x,y)\leq \rho\}\subseteq \{x\in X : J(x,y)+\lambda\Phi(x)\leq \rho+\lambda r\}\ .$$
By assumption, the set on the right-hand side is compact. Hence, the set $\{x\in S : J(x,y)\leq \rho\}$, being closed, is compact too. Summarizing:
for each $y\in Y$, the restriction of the function $J(\cdot, y)$ to $S$ is lower semicontinuous, inf-compact and has a unique global minimum. So, 
$J_{|S\times Y}$ satisfies
the hypoteses of Theorem 1.A and hence $(2.2)$ follows.\hfill $\bigtriangleup$\par
\medskip
REMARK 2.1. - From the above proof, it follows that, when $X$ is Hausdorff and each sequentially compact subset of $X$ is compact, Theorem 2.1 is still valid if we replace
``lower semicontinuous", ``inf-compact", ``closed" with ``sequentially lower semicontinuous", ``sequentially inf-compact", ``sequentially closed", respectively.\par
\medskip
We now draw a series of consequences from Theorem 2.1.\par
\medskip
COROLLARY 2.1. - {\it In addition to the assumptions of Theorem 2.1, suppose that $\beta=\sup_X\Phi$ and that $\Phi$ has no global maximum. Moreover,
suppose that the function $J(x,\cdot)$ is upper semicontinuous for all $x\in X$ and $J(x_0,\cdot)$ is sup-compact for some $x_0\in X$.\par
Then, one has
$$\sup_Y\inf_XJ=\inf_X\sup_YJ\ .$$}\par
\smallskip
PROOF. Clearly, the assumptions imply that
$$X=\bigcup_{\alpha<r<\beta}\Phi^{-1}(]-\infty,r])\ .$$
Since the family $\{\Phi^{-1}(]-\infty,r])\}_{r\in ]\alpha,\beta[}$ is filtering with respect to inclusion,
the conclusion follows from a joint application of Theorem 2.1 and Proposition 2.1 of [3].\hfill $\bigtriangleup$\par
\medskip
Another corollary of Theorem 2.1 is as follows:\par
\medskip
COROLLARY 2.2. - {\it Besides $(a)$, assume that $X$ is a convex set in a real vector space and that:\par
\noindent
$(b_1)$\hskip 5pt $\Phi$ is lower semicontinuous and strictly convex\ ;\par
\noindent
$(b_2)$\hskip 5pt  for each $\lambda>0$ and each $y\in Y$, the function $J(\cdot,y)$ is convex, lower semicontinuous and has no global minima, and
the function $J(\cdot,y)+\lambda\Phi(\cdot)$ is inf-compact.\par
Then,  for each $r\in ]\inf_X\Phi,\sup_X\Phi[$ and for each closed set $S\subseteq X$  satisfying
$$\Phi^{-1}(r)\subseteq S\subseteq \Phi^{-1}(]-\infty,r])\ ,$$
one has
$$\sup_Y\inf_SJ=\inf_S\sup_YJ\ .$$}\par
\smallskip
PROOF. We apply Theorem 2.1 taking $a=0$ and $b=+\infty$.
So, we have
 $$\alpha=\inf_X\Phi$$
as well as
$$\beta=\sup_X\Phi$$
since $M_{0,y}=\emptyset$ for all $y\in Y$. By strict convexity, the function $J(\cdot,y)+\lambda\Phi(\cdot)$ has a unique global minimum for
all $y\in Y$, $\lambda>0$. So, each assumption of Theorem 2.1 is satisfied and the conclusion follows.\hfill $\bigtriangleup$\par
\medskip
REMARK 2.2. - We stress that, in Corollary 2.2, no relation is required between the considered topology on $X$ and the algebraic structure
of the vector space which contains it.\par
\medskip
REMARK 2.3. - In the setting of Corollary 2.2, although $J$ is convex in $X$, we cannot apply the classical Fan-Sion theorem when $S$ is not convex.\par
\medskip
If $E, F$ are Banach spaces and $A\subseteq E$, a function $\psi:A\to F$ is said to be $C^1$ if it is the restriction to $A$ of a $C^1$ function on an open convex set
containing $A$.\par
\medskip
A further remarkable corollary of Theorem 2.1 is as follows:\par
\medskip
COROLLARY 2.3. - {\it Besides $(a)$, assume that $X$ is a closed and convex set in a reflexive real Banach space $E$ and that:\par
\noindent 
$(c_1)$\hskip 5pt $\Phi$ is of class $C^1$ and there is $\nu>0$ such that
$$(\Phi'(x)-\Phi'(u))(x-u)\geq\nu\|x-u\|^2$$
for all $x, u\in X$\ ;\par
\noindent
$(c_2)$\hskip 5pt  for each $y\in Y$, the function $J(\cdot,y)$ is $C^1$, sequentially weakly lower semicontinuous and $J'_x(\cdot,y)$ is Lipschitzian
with constant $L$ (independent of $y$)\ ;\par
\noindent
$(c_3)$\hskip 5pt  $\inf_{y\in Y}\inf_{M_{{{{L}\over {\nu}}},y}}\Phi>\inf_X\Phi$\ .\par
Then, for each $r\in \left ]\inf_X\Phi,\inf_{y\in Y}\inf_{M_{{L\over \nu},y}}\Phi\right [$ and for each sequentially weakly closed set $S\subseteq X$  satisfying
$$\Phi^{-1}(r)\subseteq S\subseteq \Phi^{-1}(]-\infty,r])\ ,$$
one has
$$\sup_Y\inf_SJ=\inf_S\sup_YJ\ .$$}\par
\smallskip
PROOF. For each $x, u\in X$, $y\in Y$, $\lambda\geq {{L}\over {\nu}}$, we have
$$(J'_x(x,y)+\lambda\Phi'(x)-J'_x(u,y)-\lambda\Phi'(u))(x-u)$$
$$\geq \lambda\nu\|x-u\|^2-\|J'_x(x,y)-J'_x(u,y)\|_{E^*}\|x-u\|\geq
(\lambda\nu-L)\|x-u\|^2\ .$$
Hence, the function $J(\cdot,y)+\lambda\Phi(\cdot)$, if $\lambda> {{L}\over {\nu}}$, 
is strictly convex and coercive when $X$ is unbounded ([4], pp. 247-249). Hence, if we consider $X$ with the relative weak topology, we can apply
Theorem 2.1 (in the sequential form pointed out in Remark 2.1)
 taking $a={{L}\over {\nu}}$ and $b=+\infty$, and the conclusion follows.\hfill $\bigtriangleup$\par
\medskip
If $E$ is a normed space, for each $r>0$, we put
$$B_r=\{x\in E : \|x\|\leq r\}\ .$$
If $A\subseteq E$, a function $f:A\to E$ is said to be sequentially weakly-strongly continuous if, for each $x\in A$ and for each sequence $\{x_k\}$ in $A$
weakly converging to $x$, the sequence $\{f(x_k)\}$ converges strongly to $f(x)$.\par
\medskip
COROLLARY 2.4. - {\it Let $E$ be a real Hilbert space and let $X=B_{\rho}$ for some $\rho>0$.
Besides $(a)$ and $(c_2)$, assume that
$$\delta:=\inf_{y\in Y}\|J'_x(0,y)\|>0\ .$$\
Then, for each $r\in \left ] 0,\min\left \{\rho, {{\delta}\over {2L}}\right \}\right [$, one has
$$\sup_Y\inf_{B_r}J=\inf_{B_r}\sup_YJ\ .$$}\par
\smallskip
PROOF. Apply Corollary 2.3, taking $\Phi(x)=\|x\|^2$. Let $y\in Y$ and $\tilde x\in M_{{{L}\over {2}},y}$, with $\|\tilde x\|<\rho$. Then, we have
$$J'_x(\tilde x,y)+L\tilde x=0\ .$$
Consequently, in view of $(c_2)$, we have
$$\|L\tilde x+J'_x(0,y)\|\leq \|L\tilde x\|\ .$$
In turn, using the Cauchy-Schwarz inequality, this readily implies that
$$\|\tilde x\|\geq {{\|J'_x(0,y)\|}\over {2L}}\geq {{\delta}\over {2L}}\ .$$
Therefore, we have the estimate
$$\inf_{y\in Y}\inf_{x\in M_{{{L}\over {2}},y}}\|x\|\geq \min\left \{\rho, {{\delta}\over {2L}}\right \}$$
and the conclusion follows from Corollary 2.3. \hfill $\bigtriangleup$
\medskip
We now apply Corollary 2.4 to a particular function $J$.\par
\medskip 
COROLLARY 2.5. - {\it Let $E, X$ be as in Corollary 2.4, let $Y\subseteq E$ be a closed bounded convex set and let $f:X\to E$ be a sequentially weakly-strongly continuous $C^1$ function whose derivative is Lipschitzian with constant $\gamma$. Moreover,
let $\eta$ be the Lipschitz constant of the function $x\to x-f(x)$, set
$$\theta:=\sup_{x\in X}\|f'(x)\|_{{\cal L}(E)}\ ,$$
$$L:=2\left (\eta+\theta+\gamma\left (\rho+\sup_{y\in Y}\|y\|\right )
\right )$$
and assume that
$$\sigma:=\inf_{y\in Y}\sup_{\|u\|=1}|\langle f'(0)(u),y\rangle-\langle f(0),u\rangle|>0\ .$$
Then, for each $r\in \left ] 0,\min\left \{\rho, {{\sigma}\over {L}}\right \}\right [$ and for each non-empty closed convex set $T\subseteq Y$, there exist
$x^*\in \partial B_r$ and $y^*\in T$ such that
$$\|x^*-f(x^*)\|^2+\|f(x)-y^*\|^2-\|x-f(x)\|^2\leq \|f(x^*)-y^*\|^2=(\hbox {\rm dist}(f(x^*),T))^2$$
for all $x\in B_r$\ .}\par
\smallskip
PROOF. Consider the function $J:X\times Y\to {\bf R}$ defined by
$$J(x,y)=\|f(x)-x\|^2-\|f(x)-y\|^2$$
for all $x\in X$, $y\in Y$. Clearly, for each $y\in Y$, $J(\cdot,y)$ is sequentially weakly lower semicontinuos and $C^1$. Moreover, one has
$$\langle J'_x(x,y),u\rangle = 2\langle x-f(x),u\rangle - 2\langle f'(x)(u),x-y\rangle$$
for all $x\in X$, $u\in E$. Fix $x, v\in X$ and $u\in E$, with $\|u\|=1$. We have
$${{1}\over {2}}|\langle J'_x(x,y)-J'_x(v,y),u\rangle|=|\langle x-f(x)-v+f(v),u\rangle-\langle f'(x)(u),x-y\rangle+\langle f'(v)(u),v-y\rangle|$$
$$\leq \eta\|x-v\|+|\langle f'(x)(u),x-v\rangle+\langle f'(x)(u)-f'(v)(u),v-y\rangle|$$
$$\leq \eta\|x-v\|+\|f'(x)(u)\|\|x-v\| + \|f'(x)(u)-f'(v)(u)\|\|v-y\|\leq \left (\eta+\theta+\gamma\left ( \rho+\sup_{y\in Y}\|y\|\right )\right )\|x-v\|\ .$$
Therefore, the function $J'(\cdot,y)$ is Lipschitzian with constant $L$. Fix $r\in \left ] 0,\min\left \{\rho, {{\sigma}\over {L}}\right \}\right [$ and  a non-empty closed convex set $T\subseteq Y$.
 Clearly
$$\inf_{y\in T}\||J'_x(0,y)\|\geq \inf_{y\in Y}\|J'_x(0,y)\|=2\sigma$$
and
$${{\inf_{y\in T}\||J'_x(0,y)\|}\over {2L}}>r\ .$$
Then, applying Corollary 2.4 to the restriction of $J$ to $B_r\times T$, we get
$$\sup_T\inf_{B_r}J=\inf_{B_r}\sup_TJ\ .$$
By the weak compactness of $B_r$ and $T$, we then infer the existence of $x^*\in B_r$ and $y^*\in T$ such that
$$J(x^*,y)\leq J(x^*,y^*)\leq J(x,y^*)$$
for all $x\in B_r$, $y\in T$ which is equivalent to the conclusion. To show that $x^*\in\partial B_r$, notice that if $\|x^*\|<r$ then we would have $J'_x(x^*,y^*)=0$ and so
$$r<{{\sigma}\over {L}}\leq {{\|J'_x(0,y^*)\|}\over {2L}}\leq {{L\|x^*\|}\over {2L}}<r\ ,$$
an absurd. \hfill $\bigtriangleup$
\medskip
From Corollary 2.5, in turn, we draw the following characterization about the existence and uniqueness of fixed points:\par
\medskip
COROLLARY 2.6. - {\it Let the assumptions of Corollary 2.5 be satisfied.\par
Then, for each  $r\in \left ] 0,\min\left \{\rho, {{\sigma}\over {L}}\right \}\right [$ such that $f(B_r)\subseteq Y$, 
 the following assertions are equivalent:\par
\noindent
$(i)$\hskip 5pt  the function $f$ has a unique fixed point in $B_r$ and this lies in $\partial B_r$\ ;\par
\noindent
$(ii)$\hskip 5pt the function $f$ has a fixed point in $\partial B_r$\ ;\par
\noindent
$(iii)$\hskip 5pt for each $x\in \partial B_r$ for which $f(x)\not\in B_r$, there exists $v\in B_r$ such that 
$$\|f(x)-x\|^2>\|f(v)-v\|^2-\|f(v)-f(x)\|^2\ .$$}\par
\smallskip
PROOF. The implications $(i)\to (ii)\to (iii)$ are obvious. So, suppose that $(iii)$ holds. Apply Corollary 2.5 taking $T=\overline {\hbox {\rm conv}}(f(B_r))$.
Let $x^*, y^*$ be as in the conclusion of Corollary 2.5.
Then, we have
$$\|f(x^*)-y^*\|=\hbox {\rm dist}(f(x^*),T)=0$$ and
$$\|x^*-f(x^*)\|^2+\|f(x)-f(x^*)\|^2-\|x-f(x)\|^2\leq 0 \eqno{(2.3)}$$
for all $x\in B_r$. Clearly, in view of $(iii)$, we have $f(x^*)\in B_r$. 
So, in particular, $(2.3)$ holds for $x=f(x^*)$ and this implies that
$$\|x^*-f(x^*)\|\leq 0$$
that is $x^*$ is a fixed point of $f$ in $\partial B_r$. Finally, if $\tilde x\in B_r$ and $\tilde x=f(\tilde x)$, from $(2.3)$ it follows that $f(\tilde x)=f(x^*)$,
 and so $\tilde x=x^*$. That is, $x^*$ is the unique fixed point of $f$ in $B_r$. \hfill $\bigtriangleup$\par
\medskip
REMARK 2.4. - It is important to notice that, when dim$(E)<\infty$, Corollaries 2.4, 2.5 and 2.6 are still valid replacing $B_r$ with any closed
set $S$ satisfying $\partial B_r\subseteq S\subseteq B_r$.\par
\bigskip
\bigskip
{\bf Acknowledgement.} The author has been supported by the Gruppo Nazionale per l'Analisi Matematica, la Probabilit\`a e 
le loro Applicazioni (GNAMPA) of the Istituto Nazionale di Alta Matematica (INdAM) and by the Universit\`a degli Studi di Catania, ``Piano della Ricerca 2016/2018 Linea di intervento 2". \par

\vfill\eject
\centerline {\bf References}\par
\bigskip
\bigskip
\noindent
[1]\hskip 5pt B. RICCERI, {\it Well-posedness of constrained minimization problems via saddle-points}, J. Global Optim., {\bf 40} (2008), 389-397.\par\smallskip\noindent
[2]\hskip 5pt B. RICCERI, {\it Multiplicity of global minima for parametrized functions}, Rend. Lincei Mat. Appl., {\bf 21} (2010), 47-57.\par\smallskip\noindent
[3]\hskip 5pt B. RICCERI, {\it On a minimax theorem: an improvement, a new proof and an overview of its applications},
Minimax Theory Appl., {\bf 2} (2017), 99-152.
\par
\smallskip
\noindent
[4]\hskip 5pt E. ZEIDLER, {\it Nonlinear functional analysis and its applications}, vol. III, Springer, 1985.\par

\bigskip
\bigskip
\bigskip
\bigskip
Department of Mathematics and Informatics\par
University of Catania\par
Viale A. Doria 6\par
95125 Catania, Italy\par
{\it e-mail address}: ricceri@dmi.unict.it

\bye